\documentclass[12pt,leqno]{article}
\usepackage{amssymb}
\usepackage[mathscr]{eucal}
\usepackage{amsmath,amssymb,latexsym,theorem,bbm}
\usepackage{color,url}

\newcommand{\CC}{\mathbb{C}}

\newcommand{\NN}{\mathbb{N}}

\newcommand{\RR}{\mathbb{R}}

\newcommand{\ZZ}{\mathbb{Z}}

\newcommand{\bone}{{\boldsymbol{1}}}

\newcommand{\cL}{{\mathcal L}}

\newcommand{\dd}{\mathrm{d}}
\newcommand{\ee}{\mathrm{e}}

\newcommand{\EE}{\operatorname{\mathbb{E}}}
\newcommand{\PP}{\operatorname{\mathbb{P}}}

\renewcommand{\mid}{\,|\,}

\renewcommand{\leq}{\leqslant}
\renewcommand{\geq}{\geqslant}

\newcommand{\Ci}{\mathrm{Ci}}
\newcommand{\Si}{\mathrm{Si}}
\newcommand{\si}{\mathrm{si}}

\numberwithin{equation}{section}

\theoremstyle{change} \theorembodyfont{\em}

\newtheorem{Thm}{Theorem.}[section]

\theorembodyfont{\rm}

\begin{document}

\begin{center}
{\bfseries\Large A series of definite integrals\\[1mm]
                 involving upper incomplete Gamma functions}
\\[1mm]
{\bfseries\Large }
\vspace*{3mm}

 {\sc\large
  M\'aty\'as $\text{Barczy}^{*,\diamond}$, Istv\'an $\text{Mez\H{o}}^{**}$ }

\end{center}

\vskip0.2cm

\noindent
 * HUN-REN--SZTE Analysis and Applications Research Group,
   Bolyai Institute, University of Szeged,
   Aradi v\'ertan\'uk tere 1, H--6720 Szeged, Hungary.

\noindent
 ** Department of Mathematical Sciences, University of Puerto
Rico at Mayagüez, Puerto Rico, 00681-9018, USA.

\noindent E-mails: barczy@math.u-szeged.hu (M. Barczy),
                  istvanmezo81@gmail.com  (I. Mez\H o).

\noindent $\diamond$ Corresponding author.

\renewcommand{\thefootnote}{}
\footnote{M\'aty\'as Barczy was supported by the project TKP2021-NVA-09.
Project no.\ TKP2021-NVA-09 has been implemented with the support
 provided by the Ministry of Culture and Innovation of Hungary from the National Research, Development and Innovation Fund,
 financed under the TKP2021-NVA funding scheme.}

\vspace*{-0.5cm}

\begin{abstract}
Using probability theory we derive an expression for the sum of a series of definite integrals involving upper incomplete Gamma functions.
In the proof, a normal variance mixture distribution with Beta mixing distributions plays a crucial role.
We also give an interesting application of our result, namely, a new summation formula for some derivatives
 of the Bessel functions of the first kind and the Struve functions with respect to the order.
\end{abstract}

\section{Introduction and results}

Deriving formulas for definite integrals (which, for example, involves Gamma functions) using probabilistic techniques has a long history.
We mention some recent works, although this is absolutely not an exhaustive list.
Simplified as well as original closed-form expressions for integrals involving the upper and lower incomplete Gamma functions
 were obtained by Labib and de Montigny \cite{LabMon}.
They considered series and parallel systems in which independent components have a probability of failure
 $\PP(\xi\leq t)$ up to time $t\geq 0$, where $\xi$ is a Gamma distributed random variable.
Among others, they calculated the mean-time-to-failure for a series system having
 $n$ independent components, where $n$ is a positive integer (see, Labib and de Montigny \cite[Lemmas 2.1 and 3.2]{LabMon}).

Gaunt \cite{Gau1} derived new integral formulas involving the Meijer $G$-function (see, e.g., Gradshteyn and Ryzhyik \cite[3.389(2)]{GraRyz}),
 which is defined as the contour integral of a fraction of some products of the values of the Gamma function.
In his proofs, Gaunt used a characterising equation for the product of $n$ independent Gamma distributed random variables, where $n$ is a positive integer
 (see, Gaunt \cite[Proposition 2.1]{Gau1}), and an associated distributional transformation, the gamma bias transformation (see Gaunt \cite[Section 3]{Gau1}).

Recently, Jurek \cite{Jur} has derived formulas for some definite integrals using characteristic functions
 of self-decomposable random variables.
He took advantage of the fact that a selfdecomposable random variable $\xi$ can be represented in distribution as an integral of $\ee^{-t}$
 with respect to an appropriate L\'evy process $(Y_t)_{t\geq 0}$, yielding a connection between the characteristic functions of $\xi$
 and that of the so-called background driving random variable $Y_1$ of $\xi$.

Very recently, Gaunt \cite{Gau2} has used normal variance mixture distributions with Beta and Pareto mixing distributions
 for calculating some definite integrals involving incomplete Gamma functions.
His method is based on deriving two expressions for the characteristic function of such a distribution.
This paper of Gaunt \cite{Gau2} motivated us in the sense that we also use normal variance mixture distributions with Beta mixing distributions
 in our forthcoming proof (but we do not work with characteristic functions).

Throughout this paper, let $\NN$, $\ZZ$, $\RR$, $\RR_+$, $\RR_{++}$ and $\CC$ denote the sets of positive integers, integers,
 real numbers, non-negative real numbers, positive real numbers and complex numbers, respectively.

Let $\Gamma : \RR\times \RR_+\to\RR$, $\Gamma(p,x):=\int_x^\infty u^{p-1}\ee^{-u}\,\dd u$, $p\in\RR$, $x\in\RR_+$,
 be the upper incomplete Gamma function (see, e.g., Gradshteyn and Ryzhik \cite[formula 8.350/2]{GraRyz}).
For $p>0$, we have that $\Gamma(p,0)=\int_0^\infty u^{p-1}\ee^{-u}\,\dd u$, which is nothing else but the value of the (usual)
 Gamma function at $p$, and it is denoted by $\Gamma(p)$.
For $p \in(-\infty,0)\setminus \{-n: n\in\NN\}$, the value $\Gamma(p)$ can be defined using the recursion
 $\Gamma(x+1) = x\Gamma(x)$, $x>0$, yielding, by induction, that
 \[
  \Gamma(p) = \int_0^\infty u^{p-1} \left(\ee^{-u}  - \sum_{j=0}^{n-1} \frac{(-u)^j}{j!}\right)\dd t
 \]
 for $p\in (-n,-n+1)$, $n\in\NN$, see, e.g., Fisher and Kili\c{c}man \cite[page 173]{FisKil}.
It is known that the usual Gamma function can be uniquely extended as an analytic function onto
 $\CC\setminus\{-n : n\in\ZZ\}$, and the elements of $\{-n : n\in\ZZ\}$ are all simple poles.
Further, due to Gauss's product representation, we have
 $\frac{1}{\Gamma(z)} = \lim_{n\to\infty} \frac{n^{-z}}{n!} z(z+1)\cdots(z+n)$, $z\in\CC$,
 and hence the function $\CC\ni z\mapsto \frac{1}{\Gamma(z)}$ has zeros $-n$, $n\in\ZZ_+$.

To formulate our main result, we still need to recall the definitions of the sine and cosine integrals:
  \[
   \Si(z):=\int_0^z\frac{\sin(t)}{t}\dd t,\qquad z>0,  \qquad \text{and}\qquad \Ci(z):= - \int_z^\infty \frac{\cos(t)}{t}\dd t, \qquad z>0,
  \]
  see, e.g., Abramowitz and Stegun \cite[formulas 5.2.1 and 5.2.27]{AS}.

\begin{Thm}\label{Thm_main}
The following assertions hold:
 \begin{itemize}
    \item[(i)]
   For all $a\in(0,1]$, $b\in(0,\infty)$ and $c\in\big[\frac{1}{a},\infty)$, we have
 \begin{align}\label{main_result}
  \begin{split}
  \sum_{n=1}^\infty \frac{1}{c^n} \int_0^\infty \ee^{ax-\frac{b}{x}} x^{n-\frac{3}{2}} \Gamma(1-n,x)\,\dd x
      & = 2\sqrt{\pi} \int_{\sqrt{1-a}}^\infty \frac{\ee^{-2 \sqrt{b} u}}{cu^2+ ac - 1}\,\dd u\\
      & = \sqrt{\frac{\pi}{b}} \EE\left( (c\xi^2 + ac -1)^{-1}\bone_{\{ \xi\geq \sqrt{1-a}\}} \right),
  \end{split}
 \end{align}
 where $\xi$ is an exponentially distributed random variable with parameter $2\sqrt{b}$.
   \item[(ii)]
 The series in \eqref{main_result} is convergent for $a\in(0,1)$ and $c\in\big[\frac{1}{a},\infty)$ and for $a=1$ and $c\in(1,\infty)$;
 and is divergent for $a=1$ and $c=1$.
   \item[(iii)]
For all $b\in(0,\infty)$ and $c\in\big(1,\infty)$, we have that
 \begin{align}\label{main_result_a_1}
  \begin{split}
 &\sum_{n=1}^\infty \frac{1}{c^n} \int_0^\infty \ee^{x-\frac{b}{x}} x^{n-\frac{3}{2}} \Gamma(1-n,x)\,\dd x
     = 2\sqrt{\pi}  \int_0^\infty \frac{\ee^{-2\sqrt{b} u}}{cu^2+c - 1}\,\dd u\\
 &\qquad  = \sqrt{\frac{\pi}{c(c-1)}}\Big[2\sin(\ell_{b,c})\Ci(\ell_{b,c})
             + \cos(\ell_{b,c})\left(\pi-2\Si(\ell_{b,c})\right)\Big],
  \end{split}
 \end{align}
 where $\ell_{b,c}:=2\sqrt{\frac{b(c-1)}{c}}$.
 \end{itemize}
\end{Thm}

The proof of Theorem \ref{Thm_main} can be found in Section \ref{Sect_proof1}.

According to our knowledge, the sums of the series in \eqref{main_result} and \eqref{main_result_a_1} are found neither in the literature
 nor using symbolic computation software.
In the proof of Theorem \ref{Thm_main}, a random variable $\sqrt{W}Z$ plays a crucial role,
 where $W$ and $Z$ are independent random variables such that $W$ has a Beta distribution and
 $Z$ has a standard normal distribution.
The distribution of $\sqrt{W}Z$ is called a normal variance mixture distribution with mixing distribution as that of $W$.
More generally, given $\mu,\theta\in\RR$ and $\sigma\in\RR_{++}$, the law of the random variable $\mu + \theta W+\sigma\sqrt{W}Z$
 is called a normal variance-mean mixture distribution with mixing distribution as that of $W$
 due to Barndorff-Nielsen et al.\ \cite{BarKenSor} (see also Hintz et al.\ \cite{HinHofLem}).
The proof of Theorem \ref{Thm_main} is based on the fact that we can derive two expressions
 for the expected value of $\exp\Big\{a\frac{(\sqrt{W}Z)^2}{2} - b\frac{2}{(\sqrt{W}Z)^2}\Big\}$ by two different arguments,
 where $a$ and $b$ are given in part (i) of Theorem \ref{Thm_main}.
In the first argument, we use the density function of $\sqrt{W}Z$, while, in the second argument, we use the tower rule
 by taking condition with respect to $\frac{W}{2}$.
Having these two expressions for the same expected value at hand, we can write the general term of the series in \eqref{main_result}
 as the expectation of a function of an exponentially distributed random variable,
 and then an application of the formula for the sum of a geometric series implies Theorem \ref{Thm_main}.

The approach presented in the proof of Theorem \ref{Thm_main} might be generalized in two ways
 in order to derive similar statements as in Theorem \ref{Thm_main}.
On the one hand, instead of taking $W$ as a Beta distributed random variable, one might try to use other distributions rather than the Beta distribution.
On the other hand, instead of normal variance mixture distribution, one might try to use the more general normal variance-mean mixture distributions
 for deriving similar statements as in Theorem \ref{Thm_main}.

We also give an interesting application of Theorem \ref{Thm_main}, namely,
 a new summation formula for some derivatives of the Bessel functions of the first kind and the Struve functions
 with respect to the order.
This application is seemingly completely unrelated to Theorem \ref{Thm_main}.
To describe our result, we introduce some additional notations.
The Bessel functions of the first kind are defined by
 \begin{align}\label{def_Bessel}
   J_\alpha(z):=\sum_{n=0}^\infty\frac{(-1)^n}{n!\Gamma(n+\alpha+1)}\left(\frac{z}{2}\right)^{2n+\alpha},
    \qquad z>0,\;\; \alpha\in\RR,
 \end{align}
 and the Struve functions are defined by
 \[
  H_\alpha(z):=\sum_{n=0}^\infty\frac{(-1)^n}{\Gamma\left(n+\frac32\right)\Gamma\left(n+\alpha+\frac32\right)}\left(\frac{z}{2}\right)^{2n+\alpha+1},
     \qquad z>0,\;\; \alpha\in\RR.
 \]
Further, for all $z>0$, the functions $\RR\ni\alpha\mapsto J_\alpha(z)$ and $\RR\ni\alpha\mapsto H_\alpha(z)$
 are continuously differentiable.
Recall also that $J_{-n}(z) = (-1)^n J_n(z)$, $z>0$, $n\in\NN$, and $H_{\frac{1}{2}-n}(z) = (-1)^{n-1} J_{n-\frac{1}{2}}(z)$, $z>0$, $n\in\NN$
 (see, e.g., Olver et al.\ \cite[formula 11.4.3]{OlvLozBoiCla}).
We point out that very few results are known for the derivatives of $J_\alpha$ and $H_\alpha$ with respect to the order $\alpha$.
Nonetheless, we can mention Brychkov and Geddes \cite{BryGed}, and it
 is interesting to note that an explicit evaluation has recently been found for such derivatives
 of the Bessel functions of the first kind, see the official blog of \textit{Wolfram Mathematica} \cite{MarichevBrychkov}.
Our result contributes to this small set of knowledge by providing a formula for
 the exponential generating function of the real sequence $J'_{n-\frac12}(b)-(-1)^nH'_{\frac12-n}(b)$, $n\in\NN$,
 where $b\in(0,\infty)$ and the derivative refers to the order of the Bessel and Struve functions (as elsewhere in the paper).

\begin{Thm}\label{Thm_seriesrep}
For all $b\in(0,\infty)$ and $c\in(0,1)$, we have
 \begin{align}\label{help_main2_a}
  \begin{split}
   \sum_{n=1}^\infty \frac{(bc)^n}{(n-1)!\,2^n} \Big(J'_{n-\frac12}(b)-(-1)^nH'_{\frac12-n}(b)\Big)
    =\frac{2c}{d}\sqrt{\frac{b}{2\pi}}\Big(\sin(bd)\Ci(bd)-\cos(bd)\Si(bd)\Big),
  \end{split}
 \end{align}
 where $d:=\sqrt{1-c}$, and the derivative refers to the order of the Bessel and Struve functions.
In particular, in the special case $b=2$, for all $c\in(0,1)$, we have that
 \begin{align}\label{help_main2_b}
     \sum_{n=1}^\infty\frac{c^n}{(n-1)!}\Big(J'_{n-\frac12}(2)-(-1)^nH'_{\frac12-n}(2)\Big)
         =\frac{2c}{d\sqrt\pi}\Big(\sin(2d)\Ci(2d)-\cos(2d)\Si(2d)\Big).
 \end{align}
\end{Thm}

Note that $H_{\frac{1}{2}-n}(z) = (-1)^{n-1} J_{n-\frac{1}{2}}(z)$, $z>0$, $n\in\NN$,
 however, it does not allow us to obtain a relation between $H'_{\frac{1}{2}-n}(b)$ and $J'_{n-\frac{1}{2}}(b)$ for $b\in(0,\infty)$.
The proof of Theorem \ref{Thm_seriesrep} can be found in Section \ref{Sect_proof2}.

\section{Proof of Theorem \ref{Thm_main}}
\label{Sect_proof1}

We divide the proof into eight steps.
In Steps 1-6, Step 7 and Step 8, we prove part (i), (ii) and (iii), respectively.

{\sl Step 1.}
Let $Z$ be a standard normally distributed random variable, and $W$ be a random variable having a Beta distribution
 with parameters $1$ and $\lambda$, where $\lambda\in\RR_{++}$, i.e., $W$ has a density function
 \[
    f_W(w) = \begin{cases}
              \lambda w^{\lambda-1} & \text{if $w\in(0,1)$,}\\
              0 & \text{if $w\notin(0,1)$.}
             \end{cases}
 \]
We suppose that $Z$ and $W$ are independent.
Let $X:=\sqrt{W}Z$.
Note that the conditional distribution of $X$ given $W$ is a normal distribution with mean $0$ and variance $W$.

{\sl Step 2.}
We determine the density function of $X$.
For all $x\in\RR$, by the tower rule for conditional expectations and Fubini's theorem, we have that
 \begin{align*}
   \PP(X<x) & = \PP(\sqrt{W}Z < x)
              = \EE\big( \PP\big(\sqrt{W}Z < x \mid W\big) \big)
              = \EE\left( \int_{-\infty}^x \frac{1}{\sqrt{2\pi W}}\ee^{-\frac{t^2}{2W}} \,\dd t \right) \\
            & = \int_{-\infty}^x \EE\left( \frac{1}{\sqrt{2\pi W}}\ee^{-\frac{t^2}{2W}} \right) \,\dd t.
 \end{align*}
Hence, the density function $f_X$ of $X$ takes the form
 \begin{align*}
   f_X(x) & = \EE\left( \frac{1}{\sqrt{2\pi W}}\ee^{-\frac{x^2}{2W}} \right)
            = \int_0^1 \frac{1}{\sqrt{2\pi w}} \ee^{-\frac{x^2}{2w}} \lambda w^{\lambda-1}\,\dd w
            = \frac{\lambda}{\sqrt{2\pi}}
             \int_{\frac{x^2}{2}}^\infty \left(\frac{x^2}{2y}\right)^{-\frac{1}{2}}
               \ee^{-y} \left(\frac{x^2}{2y}\right)^{\lambda-1}
               \frac{x^2}{2y^2} \,\dd y\\
          & = \frac{\lambda}{\sqrt{2\pi}}
             \left(\frac{x^2}{2}\right)^{\lambda-\frac{1}{2}}
              \int_{\frac{x^2}{2}}^\infty y^{\frac{1}{2} - \lambda -1} \ee^{-y}\,\dd y
           = \frac{\lambda}{\sqrt{2\pi}}
              \left(\frac{x^2}{2}\right)^{\lambda-\frac{1}{2}}
             \Gamma\left(\frac{1}{2} - \lambda, \frac{x^2}{2}\right),\qquad x\in\RR,
 \end{align*}
 where at the third step we used the substitution $y:=x^2/(2w)$,
 i.e.,
 \begin{align}\label{help_density_X}
   f_X(x) = \frac{\lambda}{\sqrt{\pi}} 2^{-\lambda} (x^2)^{\lambda - \frac{1}{2}}
             \Gamma\left(\frac{1}{2} - \lambda, \frac{x^2}{2}\right),
             \qquad x\in\RR.
 \end{align}
Formula \eqref{help_density_X} can also be found in Gaunt \cite[formula (2.9)]{Gau2}.

{\sl Step 3.}
We derive an expression for $\EE\Big(\ee^{a\frac{X^2}{2} - b\frac{2}{X^2}}\Big)$ using the density function $f_X$ of $X$.
By \eqref{help_density_X} and the substitution $u:=x^2/2$, we get that
 \begin{align}\label{help_int2}
   \begin{split}
   \EE\Big(\ee^{a\frac{X^2}{2} - b\frac{2}{X^2}}\Big)
    & = \int_{-\infty}^\infty \ee^{a\frac{x^2}{2} - b\frac{2}{x^2}} f_X(x) \,\dd x
      = \int_{-\infty}^\infty \ee^{a\frac{x^2}{2} - b\frac{2}{x^2}}
       \frac{\lambda}{\sqrt{\pi}} 2^{-\lambda} (x^2)^{\lambda - \frac{1}{2}}
             \Gamma\left(\frac{1}{2} - \lambda, \frac{x^2}{2}\right)\,\dd x \\
    & = \frac{\lambda 2^{1-\lambda}}{\sqrt{\pi}}
        \int_0^\infty \ee^{a\frac{x^2}{2} - b\frac{2}{x^2}}
                      (x^2)^{\lambda - \frac{1}{2}}
                       \Gamma\left(\frac{1}{2} - \lambda, \frac{x^2}{2}\right)\,\dd x \\
    & = \frac{\lambda}{\sqrt{\pi}}
        \int_0^\infty \ee^{au - b\frac{1}{u}}
          u^{\lambda - 1}
         \Gamma\left(\frac{1}{2} - \lambda, u\right)
         \,\dd u.
 \end{split}
 \end{align}

{\sl Step 4.}
We derive another expression $\EE\Big(\ee^{a\frac{X^2}{2} - b\frac{2}{X^2}}\Big)$ using the tower rule by taking condition with respect to $\frac{W}{2}$.
Namely, we obtain
 \begin{align}\label{help_cond}
 \EE\Big(\ee^{a\frac{X^2}{2} - b\frac{2}{X^2}}\Big)
  = \EE\Big(\ee^{a\frac{WZ^2}{2} - b\frac{2}{WZ^2}}\Big)
  = \EE\left(\EE\left( \ee^{a\frac{WZ^2}{2} - b\frac{2}{WZ^2}} \;\Big\vert \; \frac{W}{2}\right) \right)
  = \EE\left(h\left(\frac{W}{2}\right)\right),
 \end{align}
 where $h:(0,\frac{1}{2})\to\RR$,
 \begin{align*}
   h(u) := \int_{-\infty}^\infty \ee^{auz^2 - \frac{b}{uz^2}} \frac{1}{\sqrt{2\pi}} \ee^{-\frac{z^2}{2}}\,\dd z
         =  \frac{2}{\sqrt{2\pi}} \int_0^\infty \ee^{ -\left(\frac{1}{2} - au \right) z^2  - \frac{b}{uz^2}}\,\dd z,
         \qquad u\in\left(0,\frac{1}{2}\right).
 \end{align*}
Recall that for all $A,B\in\RR_{++}$, we have
 \begin{align}\label{help_Grad_Rhyz}
   \int_0^\infty \ee^{-Ax^2 - \frac{B}{x^2}} \,\dd x
       = \frac{1}{2}\sqrt{\frac{\pi}{A}}\,\ee^{-2\sqrt{AB}},
 \end{align}
 see, e.g., Gradshteyn and Ryzhik \cite[formula 3.325]{GraRyz}
 or Boros and Moll \cite[Example 8.4.1]{BorMol}.
For historical fidelity, we note that the evaluation of the integral in \eqref{help_Grad_Rhyz}
 is originated to Laplace.
Using \eqref{help_Grad_Rhyz} with $A:=\frac{1}{2}-au$ and $B:=\frac{b}{u}$, we have
 \[
    h(u) = \frac{1}{\sqrt{1-2au}} \exp\left(-2\sqrt{\left(\frac{1}{2}-au\right)\frac{b}{u}}\right),
          \qquad u\in\left(0,\frac{1}{2}\right).
 \]
Consequently, using also \eqref{help_cond}, we get that
 \begin{align*}
  \EE\Big(\ee^{a\frac{X^2}{2} - b\frac{2}{X^2}}\Big)
   & = \EE \left( \frac{1}{\sqrt{1-2a\frac{W}{2}}}  \exp\left(-2\sqrt{\left(\frac{1}{2}-a\frac{W}{2}\right)\frac{2b}{W}}\right)  \right) \\
   & = \EE \left( \frac{1}{\sqrt{1-aW}}  \exp\left(-2\sqrt{b}\sqrt{\frac{1-aW}{W}} \right) \right) \\
   & = \int_0^1 \frac{1}{\sqrt{1-aw}}  \exp\left(-2\sqrt{b}\sqrt{\frac{1-aw}{w}}\right)  \lambda w^{\lambda-1}\,\dd w \\
   & =  \lambda \int_0^1 \frac{w^{\lambda-\frac{3}{2}}}{\sqrt{\frac{1-aw}{w}}}  \ee^{-2\sqrt{b}\sqrt{\frac{1-aw}{w}}}  \,\dd w.
 \end{align*}
Then, by the substitution $t:=\sqrt{\frac{1-aw}{w}}$, i.e., $w=(a+t^2)^{-1}$, we have that
 \begin{align}\label{help_int1}
  \EE\Big(\ee^{a\frac{X^2}{2} - b\frac{2}{X^2}}\Big)
    = 2\lambda \int_{\sqrt{1-a}}^\infty (a+t^2)^{-\lambda-\frac{1}{2}}\ee^{-2\sqrt{b}t}\,\dd t 
    = \frac{\lambda}{\sqrt{b}} \EE\left( (a+\xi^2)^{-\lambda-\frac{1}{2}}  \bone_{\{\xi \geq \sqrt{1-a}\}}\right),
 \end{align}
 where $\xi$ is a random variable having exponential distribution with parameter $2\sqrt{b}$.
We note that, in case of $a=1$, the integral $\int_0^\infty (1+t^2)^{-\lambda-\frac{1}{2}}\ee^{-2\sqrt{b}t}\,\dd t$ in \eqref{help_int1}
 can be written in terms of Gamma function, Struve function and Bessel function of the second kind (also called Neumann function),
 see, e.g., Gradshteyn and Ryzhik \cite[formula 3.387/7]{GraRyz}.
However, we do not need such a formula in this proof.

{\sl Step 5.}
By \eqref{help_int2} and \eqref{help_int1}, we have
 \begin{align}\label{help_int3}
  \int_0^\infty \ee^{au - \frac{b}{u}}
          u^{\lambda - 1}
         \Gamma\left(\frac{1}{2} - \lambda, u\right)
         \,\dd u
  = \sqrt{\frac{\pi}{b}} \EE\left( (a +\xi^2)^{-\lambda-\frac{1}{2}} \bone_{\{\xi \geq \sqrt{1-a}\}} \right),\qquad \lambda\in\RR_{++},
 \end{align}
 where $\xi$ is a random variable having exponential distribution with parameter $2\sqrt{b}$.

{\sl Step 6.}
By \eqref{help_int3} with $\lambda:=n-\frac{1}{2}$, $n\in\NN$, and Tonelli's theorem, we get, for all $c\in\big[\frac{1}{a},\infty\big)$,
 \begin{align*}
   &\sum_{n=1}^\infty \frac{1}{c^n} \int_0^\infty \ee^{ax-\frac{b}{x}} x^{n-\frac{3}{2}} \Gamma(1-n,x)\,\dd x
      = \sqrt{\frac{\pi}{b}} \sum_{n=1}^\infty \frac{1}{c^n}  \EE\big( (a+\xi^2)^{-n}  \bone_{\{\xi \geq \sqrt{1-a}\}} \big)\\
   &\qquad = \sqrt{\frac{\pi}{b}} \EE\left( \sum_{n=1}^\infty (c(a+\xi^2))^{-n} \bone_{\{\xi \geq \sqrt{1-a}\}} \right)\\
   &\qquad = \sqrt{\frac{\pi}{b}} \EE\left( \frac{(c(a+\xi^2))^{-1}}{1 - (c(a+\xi^2))^{-1} }  \bone_{\{\xi \geq \sqrt{1-a}\}} \right)
     = \sqrt{\frac{\pi}{b}} \EE\Big( (c\xi^2+ac-1)^{-1}  \bone_{\{\xi \geq \sqrt{1-a}\}} \Big) \\
   &\qquad = 2\sqrt{\pi} \int_{\sqrt{1-a}}^\infty \frac{\ee^{-2\sqrt{b}x}}{cx^2+ac-1}\,\dd x,
 \end{align*}
 where, at the third equality, we used that $(c(a+\xi^2))^{-1}\in(0,1)$ almost surely (since $c\in\big[\frac{1}{a},\infty)$ and $\xi^2\in\RR_{++}$ almost surely).
This yields \eqref{main_result}.

{\sl Step 7.}
If $a\in(0,1)$ and $c=\frac{1}{a}$, then
 \begin{align*}
    \int_{\sqrt{1-a}}^\infty \frac{\ee^{-2\sqrt{b}x}}{cx^2+ac-1}\,\dd x
     = a\int_{\sqrt{1-a}}^\infty x^{-2}\ee^{-2\sqrt{b}x}\,\dd x
     \leq \frac{a}{1-a}\int_{\sqrt{1-a}}^\infty \ee^{-2\sqrt{b}x}\,\dd x
     \leq \frac{a}{1-a}\cdot\frac{1}{2\sqrt{b}}
     <\infty.
 \end{align*}
If $a\in(0,1]$ and $c\in(\frac{1}{a},\infty)$, then
 \begin{align*}
   \int_{\sqrt{1-a}}^\infty \frac{\ee^{-2\sqrt{b}x}}{cx^2+ac-1}\,\dd x
     \leq  \frac{1}{ac-1}\int_{\sqrt{1-a}}^\infty \ee^{-2\sqrt{b}x}\,\dd x
      \leq \frac{1}{2\sqrt{b}(ac-1)}<\infty. 
 \end{align*}
If $a=1$ and $c=1$, then
 \begin{align*}
    \int_{\sqrt{1-a}}^\infty \frac{\ee^{-2\sqrt{b}x}}{cx^2+ac-1}\,\dd x
     = \int_0^\infty x^{-2}\ee^{-2\sqrt{b}x}\,\dd x
     \geq \int_0^1 x^{-2}\ee^{-2\sqrt{b}x}\,\dd x
     \geq \ee^{-2\sqrt{b}} \int_0^1 x^{-2}\,\dd x
     =\infty.
 \end{align*}
Consequently, we have part (ii) of the assertion.

{\sl Step 8.}
The first equality in \eqref{main_result_a_1} follows readily from \eqref{main_result} by choosing $a=1$.
Now we turn to prove the second equality in \eqref{main_result_a_1}.
Using formula 3.354/1 with $\mu:=2\sqrt{b}$ and $\beta:=\sqrt{1-\frac{1}{c}}$ and formula 8.230/1 in Gradshteyn and Ryzhyik \cite[pages 341 and 886]{GraRyz},
 we get that
 \begin{align*}
    &2\sqrt{\pi}\int_0^\infty\frac{\ee^{-2\sqrt{b}u}}{ cu^2+c-1}\,\dd u
      = \frac{2\sqrt{\pi}}{c} \int_0^\infty\frac{\ee^{-2\sqrt{b}u}}{u^2+1-\frac{1}{c}}\,\dd u\\
    & = \frac{2\sqrt{\pi}}{c} \frac{1}{\sqrt{1-\frac{1}{c}}}
          \left[\sin\left(2\sqrt{b} \sqrt{1-\frac{1}{c}} \right)\Ci\left(2\sqrt{b} \sqrt{1-\frac{1}{c}} \right)
                -\cos\left(2\sqrt{b} \sqrt{1-\frac{1}{c}} \right)\si\left(2\sqrt{b} \sqrt{1-\frac{1}{c}} \right)\right] \\
     & = \frac{2\sqrt{\pi}}{\sqrt{c(c-1)}}
          \left[\sin(\ell_{b,c})\Ci(\ell_{b,c})
                -\cos(\ell_{b,c}) \left( \Si(\ell_{b,c})-\frac{\pi}{2}\right) \right]\\
     & = \frac{\sqrt{\pi}}{\sqrt{c(c-1)}} \Big[ 2\sin(\ell_{b,c})\Ci(\ell_{b,c})
                +\cos(\ell_{b,c}) \left( \pi - 2\Si(\ell_{b,c})\right) \Big]
 \end{align*}
 for $b\in(0,\infty)$ and $c\in(1,\infty)$, where $\si(z) := -\int_z^\infty \frac{\sin(t)}{t}\,\dd t$, $z>0$,
 and recall that $\ell_{b,c}=2\sqrt{\frac{b(c-1)}{c}}$.
Hence we get part (iii) of the assertion.

\section{Proof of Theorem \ref{Thm_seriesrep}}\label{Sect_proof2}

We divide the proof into three steps.
In all what follows, let $b\in(0,\infty)$ be arbitrarily fixed.

{\sl Step 1.}
By part (iii) of Theorem \ref{Thm_main}, after replacing $c$ with $1/c$ for the sake of simplicity,
 we have that
 \begin{align}\label{inter1}
   2\sqrt{\pi}\int_0^\infty\frac{\ee^{-2\sqrt{b}u}}{\frac{u^2}{c}+\frac1c-1}\,\dd u
    = \frac{c\sqrt\pi}{d}\Big(2\sin\big(2\sqrt{b}d\big)\Ci\big(2\sqrt{b}d\big)+\cos\big(2\sqrt{b}d\big)(\pi-2\Si\big(2\sqrt{b}d)\big)\Big)
 \end{align}
 for all $c\in(0,1)$, where recall that $d=\sqrt{1-c}$.

{\sl Step 2.}
We prove that
 \begin{align}\label{inter2}
  \begin{split}
   &\sum_{n=1}^\infty c^n\int_0^\infty \ee^{x-\frac{b}{x}}x^{n-\frac32}\Gamma(1-n,x)\,\dd x \\
   &\qquad = \pi{b^{-\frac{1}{4}}} \sum_{n=1}^\infty\frac{(\sqrt{b}c)^n}{(n-1)!}
              \Big( J'_{n-\frac12}(2\sqrt{b})- (-1)^n H'_{\frac12-n}(2\sqrt{b})- (-1)^n\,\pi J_{\frac12-n}(2\sqrt{b}) \Big)
  \end{split}
 \end{align}
 for all $c\in(0,1)$.
Recall that
 \[
 \Gamma(p,x)=\ee^{-x}x^p\int_0^\infty \ee^{-xu}(1+u)^{p-1}\,\dd u, \qquad p\in\RR, \;\; x>0,
 \]
 see, e.g., Kanemitsu and Tsukada \cite[formulas (3.55) and (3.58)]{KanTsu}.
Then, for each $n\in\NN$, the integral on the left-hand side of \eqref{inter2} can be written in the form
\begin{align}\label{inter3}
 \begin{split}
 \int_0^\infty \ee^{x-\frac{b}{x}}x^{n-\frac32}\Gamma(1-n,x)\,\dd x
  &=\int_0^\infty \left( \ee^{-\frac{b}{x}}x^{-\frac12}\int_0^\infty \ee^{-xu}(1+u)^{-n}\dd u\right) \dd x\\
  &=\int_0^\infty \left( (1+u)^{-n}\int_0^\infty \ee^{-\frac{b}{x}-xu}x^{-\frac12}\dd x\right)\dd u.
 \end{split}
 \end{align}
The inner integral on the right hand side of \eqref{inter3} can be transformed via the identity
\begin{equation}\label{formK}
 \int_0^\infty \ee^{-\frac{b}{x}-ux}x^A\,\dd x =  2 \left(\frac{b}{u}\right)^{\frac{A}{2}+\frac{1}{2}} K_{A+1}\left(2 \sqrt{bu}\right) 
\end{equation}
 for $u>0$ and $A\in \RR\setminus\ZZ$, see, e.g., formula 3.471/9 with $\beta:=b$, $\gamma:=u$ and $\nu:=A+1$ in Gradshteyn and Ryzhyik \cite{GraRyz}.
Here $K_\alpha$ is the modified Bessel function of the second kind given by
 \[
   K_\alpha(x):=\frac{\pi}{2}\cdot\frac{I_{-\alpha}(x)-I_\alpha(x)}{\sin(\pi\alpha)},
     \qquad x>0, \quad \alpha\in\RR\setminus \ZZ,
 \]
 where $I_\alpha$ is the modified Bessel function of the first kind defined by
 \[
 I_\alpha(x) := \sum_{m=0}^\infty \frac{1}{m! \Gamma(m+\alpha+1)} \left(\frac{x}{2}\right)^{2m+\alpha},
   \qquad x>0,\quad \alpha\in\RR,
 \]
 see, e.g., Lebedev \cite[formulas (5.7.1) and (5.7.2)]{Leb}.
In the case when $\alpha= 1/2$ (this happens when $A=-1/2$ in \eqref{formK}),
 we have the following well-known expression for the modified Bessel function of the second kind
 \[
  K_{\frac12}(x)=\sqrt{\frac{\pi }{2}}\frac{\ee^{-x}}{\sqrt{x}}, \qquad x>0,
  \]
  see, e.g., Abramowitz and Stegun \cite[10.2.17]{AS}.
Using it, \eqref{inter3} and \eqref{formK}, for each $n\in\NN$, we get that
 \begin{align}\label{help_int7}
  \int_0^\infty \ee^{x-\frac{b}{x}}x^{n-\frac32}\Gamma(1-n,x)\,\dd x
    =\sqrt\pi\int_0^\infty(1+u)^{-n}\frac{\ee^{-2 \sqrt{b}\sqrt u}}{\sqrt u}\dd u
    = 2\sqrt\pi\int_0^\infty(1+t^2)^{-n}\ee^{-2\sqrt{b}t}\dd t,
 \end{align}
 where the last equality follows by a substitution $\sqrt u=t$.
We calculate the integral at the right-hand side of \eqref{help_int7} in a roundabout way.
Namely, we derive formulae for $\int_0^\infty(1+t^2)^{-s}\ee^{-2\sqrt{b}t}\dd t$, where
 $s\in\RR_{++}\setminus (\frac{1}{2}+\ZZ)$ and $s\in\RR_{++}\setminus (\ZZ\cup (\frac{1}{2}+\ZZ))$, respectively,
 and then we take the limit of the later formula as $s\to n$, where $n\in\NN$,
 in order to get a formula $\int_0^\infty(1+t^2)^{-n}\ee^{-2\sqrt{b}t}\dd t$.
The reason for this roundabout can be explained as follows.
We cannot address any reference for a correct expression for the integral $\int_0^\infty(1+t^2)^{-n}\ee^{-2\sqrt{b}t}\dd t$, where $n\in\NN$,
 however, we can cite a reference for a (correct) expression for the integral
 $\int_0^\infty(1+t^2)^{-s}\ee^{-2\sqrt{b}t}\dd t$, where $s\in\RR_{++}\setminus (\frac{1}{2}+\ZZ)$.
Using it, we can derive another formula for the integral in question with $s\in\RR_{++}\setminus (\ZZ\cup (\frac{1}{2}+\ZZ))$,
 and then, using $\cL$'Hospital's rule, we can calculate its limit as $s\to n$, where $n\in\NN$.

To start the above described procedure, note that, similarly to \eqref{help_int7}, we can get
 \begin{align}\label{help_int8}
  \int_0^\infty \ee^{x-\frac{b}{x}}x^{s-\frac32}\Gamma(1-s,x)\,\dd x
    = 2\sqrt\pi\int_0^\infty(1+t^2)^{-s}\ee^{-2\sqrt{b}t}\dd t,
    \qquad s\in\RR_{++}.
 \end{align}
For any $s\in\RR_{++}\setminus (\frac{1}{2}+\ZZ)$, using formula 3.387/(7) in Gradshteyn and Ryzhyik \cite{GraRyz}
 with the choices $u:=1$, $\nu:=1-s$ and $\mu:=2\sqrt{b}$, we have that
 \begin{align*}
   &2\sqrt{\pi}\int_0^\infty(1+t^2)^{-s}\ee^{-2\sqrt{b}t}\dd t \\
   &\qquad  = \pi \Gamma(1-s) b^{\frac{s}{2} - \frac{1}{4}} \left(  H_{\frac{1}{2}-s}(2\sqrt{b})
                    - \frac{\cos\left( \pi \left( \frac{1}{2}-s \right) \right)}{\sin\left(\pi \left( \frac{1}{2}-s \right) \right)}J_{\frac{1}{2}-s}(2\sqrt{b})
                               +\frac{1}{\sin\left(\pi \left( \frac{1}{2}-s \right) \right)}J_{s-\frac{1}{2}}(2\sqrt{b}) \right)\\
   &\qquad = \pi \Gamma(1-s) b^{\frac{s}{2} - \frac{1}{4}}
                      \left(  H_{\frac{1}{2}-s}(2\sqrt{b}) - \frac{\sin(\pi s)}{\cos(\pi s)}J_{\frac{1}{2}-s}(2\sqrt{b})
                                     +\frac{1}{\cos(\pi s)}J_{s-\frac{1}{2}}(2\sqrt{b}) \right).
 \end{align*}
Consequently, for all $s\in\RR_{++}\setminus (\ZZ\cup (\frac{1}{2}+\ZZ))$, we get that
 \begin{align*}
   &2\sqrt{\pi}\int_0^\infty(1+t^2)^{-s}\ee^{-2\sqrt{b}t}\dd t\\
   & = \pi \Gamma(1-s) b^{\frac{s}{2} - \frac{1}{4}} \sin(\pi s) \left( \frac{1}{\sin(\pi s)} H_{\frac{1}{2}-s}(2\sqrt{b})
               - \frac{1}{\cos(\pi s)}J_{\frac{1}{2}-s}(2\sqrt{b})
                                     +\frac{1}{\cos(\pi s)\sin(\pi s)}J_{s-\frac{1}{2}}(2\sqrt{b})  \right)\\
   & = \pi \Gamma(1-s) b^{\frac{s}{2} - \frac{1}{4}} \sin(\pi s) \left( \csc(\pi s) H_{\frac{1}{2}-s}(2\sqrt{b}) - \sec(\pi s) J_{\frac{1}{2}-s}(2\sqrt{b})
                                     + 2\csc(2\pi s) J_{s-\frac{1}{2}}(2\sqrt{b})  \right),
 \end{align*}
 where $\csc$ and $\sec$ denote the cosecant and secant functions, respectively,
 defined by $\csc(\theta):=1/\sin(\theta)$, $\theta\in\RR\setminus (\pi\ZZ)$, and
 $\sec(\theta) := 1/\cos(\theta)$, $\theta\in \RR\setminus(\frac{\pi}{2} + \pi\ZZ)$, respectively.
Using Euler's reflection formula $\Gamma(z)\Gamma(1-z) = \frac{\pi}{\sin(\pi z)}$, $z\in\CC\setminus \ZZ$, we have that
 \begin{align*}
   &2\sqrt{\pi}\int_0^\infty(1+t^2)^{-s}\ee^{-2\sqrt{b}t}\dd t \\
   &\qquad =  \frac{\pi ^2 b^{\frac{s}{2} - \frac{1}{4}}
           \left(\csc (\pi  s) H_{\frac{1}{2}-s}(2\sqrt{b}) +2 \csc (2 \pi  s) J_{s-\frac{1}{2}}(2\sqrt{b})-\sec (\pi  s) J_{\frac{1}{2}-s}(2\sqrt{b})\right)}{\Gamma (s)}
 \end{align*}
 for all $s\in\RR_{++}\setminus (\ZZ\cup (\frac{1}{2}+\ZZ))$.
Taking into account \eqref{help_int8}, it implies that
 \begin{align}\label{help_int4}
 \begin{split}
  &\int_0^\infty \ee^{x-\frac{b}{x}}x^{s-\frac32}\Gamma(1-s,x)\,\dd x\\
  &\qquad = \frac{\pi^2 b^{\frac{s}{2} - \frac{1}{4}}
      \left(\csc (\pi  s) H_{\frac{1}{2}-s}(2\sqrt{b}) + 2 \csc (2 \pi  s) J_{s-\frac{1}{2}}(2\sqrt{b})-\sec (\pi  s) J_{\frac{1}{2}-s}(2\sqrt{b})\right)}{\Gamma (s)}
  \end{split}
 \end{align}
 for all $s\in \RR_{++}\setminus \Big(\ZZ\cup \Big(\frac{1}{2}+\ZZ\Big)\Big)$.

Given an $n\in\NN$, we take the limit of both sides of \eqref{help_int4} as $s\to n$.
In case of the left-hand side of \eqref{help_int4}, using \eqref{help_int7}, \eqref{help_int8},
 the inequality $(1+t^2)^{-s}\ee^{-2\sqrt{b}t} \leq \ee^{-2\sqrt{b}t}$, $t\in\RR_+$, $s>0$,
 and that $\int_0^\infty \ee^{-2\sqrt{b}t}\,\dd t = \frac{1}{2\sqrt{b}}  <\infty$, the dominated convergence theorem yields that
 \begin{align*}
  \lim_{s\to n} \int_0^\infty \ee^{x-\frac{b}{x}}x^{s-\frac32}\Gamma(1-s,x)\,\dd x
    & = 2\sqrt\pi \lim_{s\to n} \int_0^\infty(1+t^2)^{-s}\ee^{-2\sqrt{b} t}\dd t \\
    & = 2\sqrt\pi \int_0^\infty(1+t^2)^{-n}\ee^{-2\sqrt{b}t}\dd t
     = \int_0^\infty \ee^{x-\frac{b}{x}}x^{n-\frac32}\Gamma(1-n,x)\,\dd x.
 \end{align*}
Next, we calculate the limit of the right hand side of \eqref{help_int4} as $s\to n$, where $n\in\NN$.
The Gamma function is continuous on $\RR_{++}$, so the limit of the denominator $\Gamma(s)$ of the fraction at the right-hand side of \eqref{help_int4}
 is $\Gamma(n)$ as $s\to n$, where $n\in\NN$.
Consequently, it is enough to take the limit of the nominator of the fraction at the right-hand side of \eqref{help_int4}
 as $s\to n$, where $n\in\NN$.
Since
 \begin{align}\label{help_int9}
  \lim_{s\to n} b^{\frac{s}{2} - \frac{1}{4}}
         = b^{\frac{n}{2} - \frac{1}{4}} \qquad \text{and}\qquad \lim_{s\to n} \sec(\pi s) J_{\frac{1}{2}-s}(2\sqrt{b}) = (-1)^n J_{\frac{1}{2}-n}(2\sqrt{b}), \qquad n\in\NN,
 \end{align}
 it is enough to calculate the limit
 \[
  \lim_{s\to n} \left( \csc (\pi  s) H_{\frac{1}{2}-s}(2\sqrt{b}) +2 \csc (2 \pi  s) J_{s-\frac{1}{2}}(2\sqrt{b}) \right),\qquad n\in\NN.
 \]
For all $s\in\RR_{++}\setminus \big( \ZZ\cup \big(\frac{1}{2}+\ZZ\big) \big)$, we have
 \begin{align}\label{help_int11}
 \begin{split}
    &\csc (\pi  s) H_{\frac{1}{2}-s}(2\sqrt{b}) +2 \csc (2 \pi  s) J_{s-\frac{1}{2}}(2\sqrt{b})\\
    &\qquad = \frac{2\frac{\csc(2\pi s)}{\csc(\pi s)}J_{s-\frac{1}{2}}(2\sqrt{b}) +  H_{\frac{1}{2}-s}(2\sqrt{b}) }{\sin(\pi s)}
              = \frac{2 C(s) J_{s-\frac{1}{2}}(2\sqrt{b}) +  H_{\frac{1}{2}-s}(2\sqrt{b}) }{\sin(\pi s)},
  \end{split}
 \end{align}
 where we used the notation
 \begin{align*}
  C(s):=\frac{\csc(2\pi s)}{\csc(\pi s)} = \frac{\sin(\pi s)}{\sin(2\pi s)},
         \qquad s\in\RR_{++}\setminus \Big( \ZZ\cup \Big(\frac{1}{2}+\ZZ\Big) \Big).
 \end{align*}
We check that we can apply the $\cL$'Lopital's rule in order to calculate the limit of the right hand side of
 \eqref{help_int11} as $s\to n$, where $n\in\NN$.
For each $n\in\NN$, we have $\lim_{s\to n} \sin(\pi s) = \sin(\pi n)=0$, and
 \begin{align*}
 \lim_{s\to n}  \Big( 2 C(s) J_{s-\frac{1}{2}}(2\sqrt{b}) +  H_{\frac{1}{2}-s}(2\sqrt{b}) \Big)
    = 2 \left(  \lim_{s\to n} C(s)\right)  J_{n-\frac{1}{2}}(2\sqrt{b}) +  H_{\frac{1}{2}-n}(2\sqrt{b}),
 \end{align*}
 where, by the $\cL$'Lopital's rule, we have
 \begin{align}\label{help_int12}
  \begin{split}
   \lim_{s\to n} C(s) = \lim_{s\to n} \frac{\sin(\pi s)}{\sin(2\pi s)}
                      = \frac{1}{2}\lim_{s\to n} \frac{\cos(\pi s)}{\cos(2\pi s)}
                      = \frac{1}{2}\frac{\cos(\pi n)}{\cos(2\pi n)}
                      = \frac{1}{2} (-1)^n, \qquad n\in\NN.
  \end{split}
 \end{align}
This yields that
 \[
    \lim_{s\to n}  \Big( 2 C(s) J_{s-\frac{1}{2}}(2\sqrt{b}) +  H_{\frac{1}{2}-s}(2\sqrt{b}) \Big)
      = (-1)^n J_{n-\frac{1}{2}}(2\sqrt{b}) + H_{\frac{1}{2}-n}(2\sqrt{b}) =0,
 \]
 where the last equality follows, e.g., from Olver et al.\ \cite[formula 11.4.3]{OlvLozBoiCla}.
Hence we can apply the $\cL$'Lopital's rule in order to calculate the limit of the right hand side of
 \eqref{help_int11} as $s\to n$, where $n\in\NN$.
For doing so, we will need an expression for the derivative of $C$.
Note that, for all $s\in\RR_{++}\setminus \big( \ZZ\cup \big(\frac{1}{2}+\ZZ\big) \big)$, we have that
 \begin{align*}
  C'(s)& = \frac{ \pi \cos(\pi s) \sin(2\pi s) - 2\pi \sin(\pi s) \cos(2\pi s)}{(\sin(2\pi s))^2} \\
       & = \frac{ 2\pi (\cos(\pi s))^2 \sin(\pi s) - 2\pi \sin(\pi s) (2 (\cos(\pi s))^2 - 1)}{(\sin(2\pi s))^2} \\
       & = \frac{ -2\pi \sin(\pi s) (\cos(\pi s))^2 +  2\pi \sin(\pi s) } { (\sin(2\pi s))^2} \\
       & = \frac{2\pi (\sin(\pi s))^3 }{4 (\sin(\pi s))^2 (\cos(\pi s))^2 }
         = \frac{\pi}{2} \frac{\sin(\pi s)}{(\cos(\pi s))^2}
         = \frac{\pi}{2} \sec(\pi s) \tan(\pi s).
 \end{align*}
Hence
 \begin{align}\label{help_int10}
    \lim_{s\to n} C'(s) = \lim_{s\to n} \frac{\pi}{2}\frac{\sin(\pi s)}{(\cos(\pi s))^2} = \frac{\pi}{2}\cdot\frac{0}{(\cos(\pi n))^2} = 0,\qquad n\in\NN.
 \end{align}
Let us introduce the notations
 \begin{align*}
   &f(s):= 2 C(s) J_{s-\frac{1}{2}}(2\sqrt{b}) +  H_{\frac{1}{2}-s}(2\sqrt{b}),\qquad s\in\RR_{++}\setminus \Big( \ZZ\cup \Big(\frac{1}{2}+\ZZ\Big) \Big),\\
   &g(s):=\sin(\pi s), \qquad s\in\RR.
 \end{align*}
Then we have
 \[
     f'(s) = 2 C'(s) J_{s-\frac{1}{2}}(2\sqrt{b}) + 2 C(s) J_{s-\frac{1}{2}}'(2\sqrt{b}) - H_{\frac{1}{2}-s}'(2\sqrt{b}),
          \qquad s\in\RR_{++}\setminus \Big( \ZZ\cup \Big(\frac{1}{2}+\ZZ\Big) \Big).
 \]
Consequently, using \eqref{help_int12}, \eqref{help_int10} and that for all $z>0$, the functions $\RR\ni\alpha\mapsto J_\alpha(z)$ and $\RR\ni\alpha\mapsto H_\alpha(z)$
 are continuously differentiable, we get
 \begin{align*}
  \lim_{s\to n} f'(s)
      & = 2 \left(\lim_{s\to n} C'(s)\right) J_{n-\frac{1}{2}}(2\sqrt{b}) + 2 \left( \lim_{s\to n} C(s)\right) J_{n-\frac{1}{2}}'(2\sqrt{b})
          - H_{\frac{1}{2}-n}'(2\sqrt{b}) \\
      & = 2 \left( \lim_{s\to n} C(s)\right) J_{n-\frac{1}{2}}'(2\sqrt{b}) - H_{\frac{1}{2}-n}'(2\sqrt{b})\\
      & = (-1)^n J_{n-\frac{1}{2}}'(2\sqrt{b}) - H_{\frac{1}{2} - n}'(2\sqrt{b}),  \qquad n\in\NN.
 \end{align*}
Therefore, since $\lim_{s\to n} g'(s)=\pi\cos(\pi n) = (-1)^n \pi$, $n\in\NN$, we obtain that
 \[
   \lim_{s\to n} \frac{f'(s)}{g'(s)}
           = \frac{(-1)^n J_{n-\frac{1}{2}}'(2\sqrt{b}) - H_{\frac{1}{2} - n}'(2\sqrt{b})}{(-1)^n \pi}
           = \frac{1}{\pi} \Big(J_{n-\frac{1}{2}}'(2\sqrt{b}) - (-1)^n H_{\frac{1}{2} - n}'(2\sqrt{b}) \Big),
           \qquad n\in\NN.
 \]
To summarize, an application of the $\cL$'Lopital's rule implies that
 \[
   \lim_{s\to n} \left( \csc (\pi  s) H_{\frac{1}{2}-s}(2\sqrt{b}) +2 \csc (2 \pi  s) J_{s-\frac{1}{2}}(2\sqrt{b}) \right)
                    =\frac{1}{\pi} \Big(J_{n-\frac{1}{2}}'(2\sqrt{b}) - (-1)^n H_{\frac{1}{2} - n}'(2\sqrt{b}) \Big),\qquad n\in\NN.
 \]
Consequently, using \eqref{help_int9}, we get that
 \begin{align*}
   &\lim_{s\to n} \frac{\pi ^2 b^{\frac{s}{2} - \frac{1}{4}} \left(\csc (\pi  s)H_{\frac{1}{2}-s}(2\sqrt{b}) +2 \csc (2 \pi  s) J_{s-\frac{1}{2}}(2\sqrt{b})-\sec (\pi  s) J_{\frac{1}{2}-s}(2\sqrt{b})\right)}{\Gamma (s)}\\
   &\qquad = \frac{\pi^2b^{\frac{n}{2} - \frac{1}{4}}
       \Big(\frac{1}{\pi} J'_{n-\frac{1}{2}}(2\sqrt{b}) - \frac{1}{\pi} (-1)^n H'_{\frac{1}{2}-n}(2\sqrt{b}) - (-1)^n J_{\frac{1}{2}-n}(2\sqrt{b}) \Big)}{\Gamma(n)} \\
   &\qquad  = \frac{\pi b^{\frac{n}{2} - \frac{1}{4}} }{(n-1)!}\Big( J'_{n-\frac{1}{2}}(2\sqrt{b}) - (-1)^n H'_{\frac{1}{2}-n}(2\sqrt{b}) - (-1)^n \pi J_{\frac{1}{2}-n}(2\sqrt{b}) \Big),
    \qquad n\in\NN.
 \end{align*}
All in all, by taking the limit of both sides of \eqref{help_int4}, we arrive at
 \begin{align}\label{help_int6}
  \begin{split}
   \int_0^\infty \ee^{x-\frac{b}{x}}x^{n-\frac32}\Gamma(1-n,x)\,\dd x
        = \frac{\pi b^{\frac{n}{2} - \frac{1}{4}} }{(n-1)!}\Big( J'_{n-\frac{1}{2}}(2\sqrt{b}) - (-1)^n H'_{\frac{1}{2}-n}(2\sqrt{b}) - (-1)^n \pi J_{\frac{1}{2}-n}(2\sqrt{b}) \Big)
  \end{split}
 \end{align}
 for each $n\in\NN$.
Multiplying both sides of \eqref{help_int6} by $c^n$ and then summing up for $n\in\NN$, we have \eqref{inter2}, as desired.

{\sl Step 3.}
We derive the statement of the theorem.
By Theorem \ref{Thm_main}, the left-hand sides of  \eqref{inter1} and \eqref{inter2} coincide, and hence
 their right-hand sides coincide as well, yielding that
\begin{align}\label{SiCi}
 \begin{split}
  \frac{ b^{\frac{1}{4}} c}{d\sqrt\pi}&\left(2\sin(2\sqrt{b}d)\Ci(2\sqrt{b}d)+\cos(2\sqrt{b}d)(\pi-2\Si(2\sqrt{b}d))\right)\\
   & = \sum_{n=1}^\infty\frac{(\sqrt{b}c)^n}{(n-1)!} \Big( J'_{n-\frac12}(2\sqrt{b})- (-1)^n H'_{\frac12-n}(2\sqrt{b})- (-1)^n\,\pi J_{\frac12-n}(2\sqrt{b}) \Big)
 \end{split}
 \end{align}
 for all $c\in(0,1)$, where recall that $d=\sqrt{1-c}$.
 Next, we check that
 \begin{equation}\label{Jcos}
  \sum_{n=1}^\infty\frac{(\sqrt{b}c)^n}{(n-1)!}(-1)^nJ_{\frac12-n}(2\sqrt{b})
      = -\frac{b^{\frac{1}{4}}c}{d}\cdot\frac{\cos(2\sqrt{b}d)}{\sqrt\pi}, \qquad c\in(0,1).
 \end{equation}
Using Prudnikov et al.\ \cite[5.7.6/(2), page 660]{PruBryMar}, for all $c\in(0,1)$,
 we have that
 \begin{align*}
    \sum_{n=1}^\infty\frac{(\sqrt{b}c)^n}{(n-1)!}(-1)^nJ_{\frac12-n}(2\sqrt{b})
       & = (-\sqrt{b}c)\sum_{n=0}^\infty\frac{(-\sqrt{b}c)^n}{n!}J_{-n-\frac12}(2\sqrt{b}) \\
       &= (-\sqrt{b}c) (2\sqrt{b})^{\frac{1}{4}} (2\sqrt{b}-2\sqrt{b}c)^{-\frac{1}{4}}
              J_{-\frac{1}{2}}\left(\sqrt{4b-4bc}\right) \\
       &= - \sqrt{b}c(1-c)^{-\frac{1}{4}} J_{-\frac{1}{2}}\left(2\sqrt{b}\sqrt{1-c}\right)\\
       &=- \sqrt{b}c d^{-\frac{1}{2}} J_{-\frac{1}{2}}\left(2\sqrt{b}d\right). 
  \end{align*}
Using the explicit formula for $J_{-\frac{1}{2}}(z)$, $z\ne 0$ (see, e.g., Gradshteyn and Ryzhyik \cite[8.464/(2)]{GraRyz}),
 it implies that
 \begin{align*}
    \sum_{n=1}^\infty\frac{(\sqrt{b}c)^n}{(n-1)!}(-1)^nJ_{\frac12-n}(2\sqrt{b})
        = -\sqrt{b} c d^{-\frac{1}{2}} \sqrt{\frac{2}{\pi2\sqrt{b}d}}
                \cos(2\sqrt{b}d)
      = -\frac{b^{\frac{1}{4}}c}{d}\cdot\frac{\cos(2\sqrt{b}d)}{\sqrt\pi}
 \end{align*}
 for all $c\in(0,1)$, as desired.

Finally, using \eqref{SiCi} and \eqref{Jcos}, for all $c\in(0,1)$, we obtain that
 \begin{align*}
  & \frac{2b^{\frac{1}{4}}c}{d\sqrt\pi}\Big(\sin(2\sqrt{b}d)\Ci(2\sqrt{b}d)-\cos(2\sqrt{b}d)\Si(2\sqrt{b}d)\Big)
       + \frac{b^{\frac{1}{4}}c\sqrt{\pi}}{d}\cos(2\sqrt{b}d)  \\
  &\qquad = \sum_{n=1}^\infty\frac{(\sqrt{b}c)^n}{(n-1)!}\Big(J'_{n-\frac12}(2\sqrt{b})-(-1)^nH'_{\frac12-n}(2\sqrt{b})\Big)
            + \frac{b^{\frac{1}{4}}c\sqrt{\pi}}{d}\cos(2\sqrt{b}d),
 \end{align*}
 yielding that
 \begin{align*}
   \sum_{n=1}^\infty\frac{(\sqrt{b}c)^n}{(n-1)!}\Big(J'_{n-\frac12}(2\sqrt{b})-(-1)^nH'_{\frac12-n}(2\sqrt{b})\Big)
   = \frac{2b^{\frac{1}{4}}c}{d\sqrt\pi}\Big(\sin(2\sqrt{b}d)\Ci(2\sqrt{b}d) - \cos(2\sqrt{b}d)\Si(2\sqrt{b}d)\Big).
 \end{align*}
Replacing $b$ by $\frac{b^2}{4}$ in the previous equality, we have \eqref{help_main2_a}.

By choosing $b=2$ in \eqref{help_main2_a}, we obtain \eqref{help_main2_b} as well.
The proof of the theorem is therefore finished.

\section*{Acknowledgements}
We acknowledge the valuable suggestions from the referee.
In particular, we are grateful for drawing our attention to the possible extension
 of Theorem \ref{Thm_main} for parameters $a\in(0,1]$ and $b\in(0,\infty)$ instead of $a=b=1$.

\section*{Statements and Declarations}

The authors declare that they have no known competing financial interests
 or personal relationships that could have appeared to influence the work presented in this paper.

\end{document}